\documentclass{amsproc}

\usepackage{amssymb}
\usepackage{bm}

\usepackage{subfiles}
\usepackage{}

\newcommand{\SPAN}{\mathop\mathrm{span}}
\newcommand{\TR}{\mathop\mathrm{tr}}

\newtheorem{theorem}{Theorem}[section]

\newtheorem{prop}[theorem]{Proposition}
\newtheorem{corollary}[theorem]{Corollary}

\theoremstyle{definition}
\newtheorem{definition}[theorem]{Definition}

\theoremstyle{remark}
\newtheorem{remark}[theorem]{Remark}

\newtheorem{question}[theorem]{Question}
\newtheorem{answer}[theorem]{Answer}

\numberwithin{equation}{section}

\begin{document}

\title{Constructions of biangular tight frames and their relationships with equiangular tight frames}


\author{Peter G. Casazza}
\email{casazzap@missouri.edu}
\address{Department of Mathematics, University
of Missouri, Columbia, MO 65211-4100}
\thanks{The authors were supported by
 NSF DMS 1307685; and  NSF ATD 1321779; ARO W911NF-16-1-0008.  In addition, the authors would like to thank Dustin Mixon for helpful discussions regarding the affine picket-fence BTFs discussed in Section~\ref{sect_steiner}.}

\author{Jameson Cahill}
\email{jamesonc@nmsu.edu}
\address{Department of Mathematical Science, New Mexico State University, Las Cruces, NM 88003}

\author{John I. Haas}
\email{haasji@missouri.edu}
\address{Department of Mathematics, University
of Missouri, Columbia, MO 65211-4100}

\author{Janet Tremain}
\email{tremainjc@missouri.edu}
\address{Department of Mathematics, University
of Missouri, Columbia, MO 65211-4100}

\subjclass[2000]{Primary }

\date{}

\begin{abstract}
We study several interesting examples of {\bf Biangular Tight Frames (BTFs)}  -  basis-like sets of unit vectors admitting exactly two distinct {\bf frame angles} (ie, pairwise absolute inner products)  - and examine their relationships with {\bf Equiangular Tight Frames (ETFs)} - basis-like systems which admit exactly one frame angle (of minimal coherence).

We develop a general framework of so-called Steiner BTFs - which includes the well-known Steiner ETFs as special cases; surprisingly, the development of this framework leads to a connection with famously open problems regarding the existence of Mersenne and Fermat primes.

In addition, we demonstrate an example of a  smooth parametrization of $6$-vector BTFs in $\mathbb R^3$, where the curve ``passes through'' an ETF; moreover, the corresponding frame angles ``deform'' smoothly with the parametrization, thereby answering two questions  about the rigidity of BTFs.

Finally, we generalize from BTFs to {\bf (chordally) biangular tight fusion frames (BTFFs)} - basis-like sets of orthogonal projections admitting exactly two distinct trace inner products - and we explain how one may think of them as generalizations of  BTFs.  In particular, we construct an interesting example of a BTFF corresponding to $16$ $2$-dimensional subspaces of $\mathbb R^4$ that ``Pl{\"u}cker embeds'' into a Steiner ETF consting of $16$ vectors in $\mathbb R^6$, which  refer to as a {\bf Pl{\"u}cker ETF.}
\end{abstract}

\maketitle


\section{Introduction}\label{sec_intro}

Research into {\it finite (unit-norm) frames} - basis-like sets of unit vectors -  has developed rapidly over the last two decades~\cite{MR836025, MR2892742}. While the discipline has found various applications~\cite{MR2105392}, possibly  the most common use involves the analysis of signals~\cite{MR1486646, MR2021601, Zauner1999, 1523643}.  
Accordingly, a common criterion for frame design is low {\it coherence}, as this may lead to improved error reduction in various signal processing applications~\cite{MR2021601, MR1984549}. 

The {\it coherence} of a frame is the maximal value among its set of pairwise absolute inner products.  Frames that achieves minimal coherence are called a {\bf Grassmannian frames}, a class including the well-studied {\it equiangular tight frames (ETFs)}~\cite{Fickus:2015aa} -  which admit exactly  one pairwise absolute inner product. Recent advances have shown that {\it biangular tight frames(BTFs)} - unit-norm frames which admit exactly two pairwise absolute inner products - sometimes form Grassmannian frames~\cite{Hoggar1982, WoottersFields1989} or can be used to construct Grassmannian frames~\cite{2015arXiv150905333B, 2016arXiv161003142C}.
Accordingly, this work is a study of the relationship between ETFs and BTFs.  In order to illuminate their relationship, we consider interesting examples; moreover, we consider general frameworks for {\bf harmonic BTFs/ETFs} and so-called {\bf Steiner BTFs/ETFs}

 In Section~\ref{sect_pre}, we recall basic facts from frame theory and fix notation.  In Section~\ref{sect_continuum}, we demonstrate a  smooth parametrization of $6$ vector BTFs in $\mathbb R^3$, where the corresponding frame angles transform smoothly with the parameter.  This smooth curve of BTFs ``passes through'' an ETF, and its existence answers two questions regarding the rigidity of BTFs.  In Section~\ref{sect_harmonic} and Section~\ref{sect_steiner}, we develop a framework based on Steiner systems and divisible difference sets, which produces numerous infinite families of so-called {\bf Steiner BTFs}, which admits the {\bf Steiner ETFs}~\cite{MR2890902} as a special case. Surprisingly, we observe a connection between one of these families and the long-standing open problem(s)~\cite{MR990017} regarding the existence of an infinitude of Mersenne primes and Fermat primes.
Finally, in Section~\ref{sect_plucker}, we generalizes from BTFs to {\bf (chordally) biangular tight fusion frames (BTFFs)} - basis-like sets of projections admitting exactly two distinct trace inner products.
We define a new type of ETF - a {\bf Pl{\"u}cker ETF} - based on the {\bf Pl{\"u}cker embedding}; as an example, we construct $16$ $2$-dimensional subspaces of $\mathbb R^4$ that {\bf Pl{\"u}cker embed} into a Steiner ETF of $16$ vectors in $\mathbb R^6$.
constructing
``Pl{\"u}cker embedding'' a {\bf (chordally) biangular tight fusion frame (BTFF)}.


\section{Preliminaries}\label{sect_pre}
Given a finite dimensional Hilbert space $\mathbb F^m$, with the usual inner product and where $\mathbb F = \mathbb R$ or $\mathbb C$, we fix $\{e_j\}_{j=1}^m$ as the canonical orthonormal basis, and we fix $I_m$ as the $m \times m$ identity matrix. A set of vectors $\mathcal F = \{f_j\}_{j =1}^n \subset \mathbb F^m$ is a {\bf (finite) frame} if $\SPAN \{f_j\}_{j=1}^n = \mathbb F^m.$  
 It is convenient to think of the frame $\mathcal F$  in terms of its {\bf synthesis matrix},
$F = [ f_1 \, f_2  \, ... \, f_n ],$
the $m \times n$ matrix with columns given by the frame vectors, so we identify the set $\mathcal F$ with the matrix $F$ as representatives of the same object.
We say that $\mathcal F$ is {\bf $\boldsymbol a$-tight} if $FF^*= a I_m$ for some $a>0$, called the {\bf tightness parameter}, where, in general, $A^*$ denotes the adjoint of a matrix $A$. We say that $\mathcal F$ is {\bf unit-norm} if each frame vector has norm $\|f_j\|=1$.

To abbreviate, a unit-norm, $a$-tight frame consisiting of $n$ vectors for $\mathbb F^m$ is called an {\bf $\boldsymbol{(n,m)}$-frame}. 
The tightness parameter for such a frame is determined as $a=\frac{n}{m}$, because
$$n= \sum_{j'=1}^m\sum_{j=1}^n |\langle e_{j'}, f_j \rangle|^2 = \sum_{j'=1}^m\sum_{j=1}^n \TR(f_j  f_j^* e_{j'}  e_{j'}^*) = a \sum_{j'=1}^m \| e_{j'} \|^2 = am,$$
which also implies that every such frame satisfies the identity 
\begin{equation}\label{eq_unt_cond}
\sum\limits_{j'=1}^n |\langle f_j, f_{j'} \rangle|^2 = \frac{n}{m} \text{ for every } j \in \{1,...,n\}.
\end{equation}

Given any unit-norm frame $\mathcal F = \{ f_j \}_{j=1}^n$, its {\bf frame angles} are the elements of the frame's {\bf angle set}, which we denote and define by
$
\Theta_{\mathcal F} : =\big\{ |\langle f_j, f_{j'} \rangle | : j \neq j' \big\}.
$
We say that $\mathcal F$ is {\bf $\boldsymbol d$-angular} if $|\Theta_{\mathcal F}| =d$ for some $d \in \mathbb N$.  
In the special case that $\mathcal F$ is $1$-angular or $2$-angular, then we say that it is {\bf equiangular} or {\bf biangular}, respectively.  
If $\mathcal F$ is $d$-angular frame with frame angles $\alpha_1, \alpha_2,...,\alpha_d$,  then $\mathcal F$ is {\bf equidistributed} if there exist positive integers $\tau_1, \tau_2, ..., \tau_d \in \mathbb N$ such that 
$$
\Big| \Big\{ j' \in \small\{1,...,n \small\} : j' \neq j, |\langle f_j, f_{j'} \rangle | = \alpha_l \Big\} \Big| = \tau_l
$$
for every $j \in \{1,2,...,n\}$ and every $l \in \{1,2,...,d\}$.  In this case, we call the positive integers $\tau_1, \tau_2,...,\tau_d$  the {\bf frame angle multiplicities} of $\mathcal F$ and note that $\sum_{j=1}^d \tau_j= n-1$.

By the lower bound of Welch \cite{Welch1974}, if $\mathcal F=\{f_j\}_{j=1}^n$ is a unit-norm frame for $\mathbb F^m$, then the the frame's {\bf coherence} - the maximal element of its angle set - obeys the inequality,
\begin{equation*}
\max\limits_{j \neq j'} |\langle f_j, f_{j'} \rangle| \geq \sqrt{\frac{n-m}{m(n-1)}},
\end{equation*}
and it is well-known \cite{ Fickus:2015aa} that $\mathcal F$ achieves this bound if and only if $\mathcal F$ is an equiangular, tight frame ({\bf ETF}).  For convenience, we call this lower bound the {\bf Welch constant} and write $W_{n,m}=\sqrt{\frac{n-m}{m(n-1)}}.$

In general, $d$-angular $(n,m)$-frames are not equidistributed, as is demonstrated by the $3$-angular {\it orthoplectic} constructions in  \cite{2015arXiv150905333B}, for example.  However, the equidistributed property always holds in the special case that  $\mathcal F$ is a $d$-angular $(n,m)$-frame with $d \leq 2$. This is trivial if $d=1$, since the frame's angle set is then a singleton.    If $\mathcal F$ is a biangular tight frame ({\bf BTF}), then straightforward substitution implies the invariance \big(with respect to the choice $j\in\{1,...,n\}$\big) of the distribution of the two squared frame angles occurring as summands in the left-hand side of Equation~(\ref{eq_unt_cond});   in other words,  $\mathcal F$ is equidistributed.  Furthermore,  solving Equation~(\ref{eq_unt_cond}) in conjunction with the identity $\tau_1 + \tau_2 = n-1$ yields the frame angle multiplicities.

\begin{prop}\label{prop_btf_multiplicities}
If $\mathcal F$ is a biangular $(n,m)$-frame with distinct frame angles $\alpha_1$ and $\alpha_2$, then the correponding frame angle multiplicites are
$$\tau_1 = \frac{n-1}{\alpha_2^2 - \alpha_1^2}\left(\alpha_2^2 - \frac{n-m}{m(n-1)}\right)
 \text{ and } 
\tau_2 = \frac{n-1}{\alpha_1^2 - \alpha_2^2}\left(\alpha_1^2 - \frac{n-m}{m(n-1)} \right).
$$
\end{prop}

We refer to \cite{MR2964005} for more information about frame theory and its applications.
\begin{remark}\label{rem_abuse}
In order to provide a concise definition of the {\it equidistributed} property above, it was convenient to define  ETFs and BTFs as distinctly different objects, as determined by the cardinalities of their frame angle sets.   However, we prefer to think of ETFs as a special instances of BTFs; that is, after allowing a slight abuse of the definition, we think of an ETF as a BTF where the two frame angles agree.  For the purpose of highlighting the relationship between these two types of frames, we will employ this ``abuse of terminology'' throughout the remainder of this work, but never without clarification.
\end{remark}


\section{A continuum of BTFs in $\mathbb R^3$}\label{sect_continuum}
In this section, we construct a smooth curve of BTFs passing through the space of real $(6,3)$-frames.  Afterwards, we remark on some of the curve's surprising properties and the implications of its existence.

Define the curve
$F: [1, \infty)  \rightarrow \mathbb R^{3 \times 6}$
defined by
$$
F(t) = 
\left[
\begin{array}{cccccc}
\frac{1}{t^2} & \frac{1}{t^2}  & 0 & 0 & \sqrt{\frac{t^4 -  1}{t^4}} & -\sqrt{\frac{t^4 -  1}{t^4}} \\
0 & 0 & \sqrt{\frac{t^4 -  1}{t^4}} & -\sqrt{\frac{t^4 -  1}{t^4}} & \frac{1}{t^2} & \frac{1}{t^2} \\
 \sqrt{\frac{t^4 -  1}{t^4}} & -\sqrt{\frac{t^4 -  1}{t^4}} & \frac{1}{t^2} & \frac{1}{t^2} & 0 & 0 \\
\end{array}
\right].
$$
The following properties of $F(t)$ follow by elementary calculations:
\begin{enumerate} 
\item[{(i)}] for each $t \in [1, \infty)$, the rows are pairwise orthogonal,
\item[{(ii)}] for each $t \in [1, \infty)$, every row has squared norm equal to $2$,
 \item[{(iii)}] for each $t \in [1, \infty)$, every column has norm $1$ and
 \item[{(iv)}] for each $t \in [1, \infty)$, the angle set is
$$
\Theta(t) = \left\{\left|\frac{2-t^4}{t^4}\right|, \frac{\sqrt{t^4 -1}}{t^4} \right\}.
$$
\end{enumerate}

By our identification of frames with their synthesis operators, $F(t)$ corresponds to a smooth curve, $\mathcal F(t)$, of frames passing through the space of real $(6,3)$-frames, where every point on the curve is a BTF.
A remarkable property of this curve is that it passes through three frames of noteworthy structural symmetry, and it converges to a fourth noteworthy frame at infinity.  

The first ocurs at $t=1$, where we obtain
$$
F\big(1\big) = 
\left[
\begin{array}{cccccc}
1 & 1  & 0 & 0 & 0 & 0 \\
0 & 0 & 0 & 0 & 1 & 1 \\
0 & 0 & 1 & 1 & 0 & 0 \\
\end{array}
\right],
$$
or two copies of the canonical basis, $\{e_1,e_2,e_3\}$.  Although $F(1)$ is  an ``uninteresting'' frame, it is nevertheless highly symmetric and it provides a convenient visualization for our curve's initial point.  We imagine that as we ``trace out'' the parameter from here, this pair of orthonormal bases  smoothly deforms to become other interesting frames, while retaining the rigid BTF structure along the way. 

At $t=\sqrt[4]{2}$, we obtain the frame
$$
F\big(\sqrt[4]{2}\big) = 
\left[
\begin{array}{cccccc}
\frac{1}{\sqrt 2} & \frac{1}{\sqrt 2}  & 0 & 0 & \frac{1}{\sqrt 2} & \frac{-1}{\sqrt 2} \\
0 & 0 & \frac{1}{\sqrt 2} & \frac{-1}{\sqrt 2} & \frac{1}{\sqrt 2} & \frac{1}{\sqrt 2} \\
\frac{1}{\sqrt 2} &\frac{-1}{\sqrt 2} & \frac{1}{\sqrt 2} & \frac{1}{\sqrt 2}  & 0 & 0 \\
\end{array}
\right],
$$
 an example of Steiner BTF, as discussed in Section~\ref{sect_steiner}. 
Furthermore, by partitioning its columns, $\{f_j\}_{j=1}^6$,  into three consecutive pairs, $\{f_1, f_2\}$, $\{f_3, f_4\}$ and $\{f_5, f_6\}$, we obtain all of the orthonormal bases for the $2$-dimensional coordinate subspaces of $\mathbb R^3$. 

Arguably, the most noteworthy point~\cite{Fickus:2015aa, BK06} on the curve occurs at $t=\sqrt[4]{\frac{5 +\sqrt{5}}{2}}$, because the angle set reduces to the singleton,
$$
\Theta\big(\scriptstyle{\sqrt[4]{\frac{5 +\sqrt{5}}{2}}} \big) = \big\{ \scriptstyle \frac{1}{\sqrt 5} \big\},
$$
meaning that $F\big(\scriptstyle{\sqrt[4]{\frac{5 +\sqrt{5}}{2}}} \big)$ is an ETF.  This particular ETF is well known~\cite{BK06, MR2964005}, as it corresponds to a selection of six non-antipodal  vertices from a regular icosahedron.  In accordance with the characterization of Welch, this frame is a minimizer of coherence among all  real $(6,3)$-frames; that is, it is a Grassmannian frame. 

Finally, we consider the asymptotics. Evaluating the limits of the coordinates as $t$ approaches infinity yields
$$
\lim_{t\rightarrow \infty} F(t)=  
\left[
\begin{array}{cccccc}
0 & 0 & 0 & 0 & 1 & -1\\
 0 & 0 & 1 & -1 & 0 & 0\\
 1 & -1 & 0 & 0 & 0 & 0
\end{array}
\right],
$$
the canonical basis, $\{e_1, e_2, e_3\}$, unioned with its antipodes, $\{-e_1, -e_2, -e_3\}$.  As with the the curve's initital point, the curve's limiting frame is not ``interesting'' but highly symmetric, thereby supplementing our visualization of how the BTFs deform as they  transit along the curve.

Thus, we have constructured a curve of highly structured frames (BTFs) whose inital point is a frame with worst possible coherence -  two copies of the canonical basis - which then passes through a frame with optimal coherence - an 
ETF -  and then asymptotically approaches another frame with worst possible coherence - the canonical basis along with its negatives.  The existence of such a curve answers two questions concerning the existence and properties of BTFs.

\subsection{Two questions answered}
As with ETFs, biangular tight frames necessitate the equidistributed property, a fairly rigid constraint.   Moreover, just as ETFs are known to be difficult to construct~\cite{MR2460526,ConwayHardinSloane1996,MR2890902,2014arXiv1408.0334H, MR2021601, MR3150919, MR1984549, SustikTroppDhillonHeath2007, Szollosi2014a, Fickus:2015aa}, attempts to construct biangular tight frames are accompanied by a similar level of difficulty~\cite{RoyScott2007, MR1383512, Neumaier1989, MR2763068,DelsarteGoethalsSeidel1975, 2014arXiv1402.3521B, Hoggar1982, WoottersFields1989, 2014arXiv1402.3521B}.  Because of these similariities, we had, informally speaking,  begun to wonder:
\begin{quote}
``How rigid are BTFs compared to ETFs?"
\end{quote}
As it turns out, the existence of our curve, $F(t)$, answers two meaningful formulations of this question. 
Due to the Welch bound, the square of an ETF's single frame angle is always rational. Based on known examples~\cite{2016arXiv161003142C,Neumaier1989, Hoggar1982, WoottersFields1989}, we had begun to suspect a similar statement for BTFs.
\begin{question}\label{quest1}
Must the squared frame angles of a BTF be rational - or, at least, nontranscendental?
\end{question}

\begin{answer}
Our curve's existence answers this in the negative. The continuum of frame angles, $\{\Theta(t) \}_{t\in[1,\infty)},$ along with elementary properties of the real line indicate the existence of uncountably many inquivalent BTFs occurring along the curve with transcendental squared frame angles.  
\end{answer}

Incidentally, this observation leads to and answers the next question.
     We have also pondered the potential cardinality of {\bf inequivalent} BTFs that might occur within a given space of $(n,m)$-frames, where two biangular $(n,m)$-frames are {\bf inequivalent} if they have different frame angle sets. 
There are several instances where two or more inequivalent BTFs coexist within a given space~\cite{2016arXiv161003142C}.  For example, there are at least three inequivalent BTFs for $\mathbb C^3$ consisting of $8$ vectors~\cite{2016arXiv161003142C}.  However, our limited knowledge of known examples~\cite{2016arXiv161003142C,Neumaier1989, Hoggar1982, WoottersFields1989 } of BTFs led us to ponder the following.
\begin{question}\label{quest2}
Within a given space of $(n,m)$ frames, must the number of inequivalent BTFs be finite - or, at least, countable?
\end{question}
\begin{answer}Again, our curve's existence answers this in the negative. The continuum of frame angles, $\{\Theta_t\}_{t\in [1,\infty)}$, implies to the existence of uncountably many inequivalent BTFs coexisting within the space of $(6,3)$-frames.\end{answer}
    Note that, since the real $(6,3)$-frames are a subset of the complex $(6,3)$-frames, we may view our curve $F(t)$ as a complex enitity.  In particular, Question~\ref{quest1} and Question~\ref{quest2} are answered in the negative for both the real and complex cases.


\section{Harmonic BTFs}\label{sect_harmonic}
Given that our thesis is to outline the peculiar connections between BTFs and ETFs, a brief discussion of their relationship when manifesting as so-called {\it harmonic frames} seems appropriate.  Besides their relevance to this paper's overall theme, we will require basic facts about harmonic ETFs and BTFs  to construct the {\it Steiner BTFs} of the next section.  Unfortunately, a proper treatment of harmonic frames requires a level of notational complexity that we prefer to avoid in this work, so this section is only intended as a summary of core definitions and results~\cite{2016arXiv161003142C, MR1984549,  XiaZhouGiannakis2005}.

In order to define {\it harmonic frames}, we assume basic familiarity with character theory (see \cite{MR1625181} for details) and recall our identification of frames with their synthesis matrices.  
Let $\mathbb G$ be an abelian group of order $n$ and let $\mathcal S$ be any subselection of $m$ distinct rows from the $n\times n$ character table of $\mathbb G$.  The resulting $m \times n$ submatrix, $H$, with columns rescaled to norm $1$, is called an {\bf $(n,m)$-harmonic frame (for $\mathbb G$ generated by $\mathcal S$)}. By basic properties of character tables~\cite{MR1625181}, it is easy to verify~\cite{MR2964005, 2016arXiv161003142C} that such an object is the sythesis matrix of an $(n,m)$-frame, $\mathcal H$, for $\mathbb C^m.$   Moreover,  every $(n,m)$-harmonic frame is {\bf flat}, meaning that its entries have constant magnitude $1 / \sqrt m$; the importance of this property becomes apparent in the next section.  Also, we remark that our choice of $\mathcal S$  - a combinatorial issue - completely determines $\mathcal H$'s frame angle set.   Assimilating these facts, assured that existence is not an issue, we offer the following simplified definition, with notational details purposefully suppressed. 
\begin{definition}
Given an abelian group $\mathbb G$ of order $n$ and a subset $\mathcal S$ of $\mathbb G$ of $m$ distinct elements, then an {\bf $(n,m)$-harmonic frame (for $\mathbb G$ generated by $\mathcal S$)} is a flat $(n,m)$-frame for $\mathbb C^m$ with a frame angle set, $\Theta(\mathbb G, \mathcal S)$, determined by $\mathbb G$ and $\mathcal S$.
\end{definition}

We remark that, under certain conditions on $\mathbb G$ or the subselection $\mathcal S$, a harmonic frame can manifest as a strictly real frame~\cite{2016arXiv161003142C}.
For more details about harmonic frames, we refer to Waldron's chapter in~\cite{MR2964005} and the references therein.

In order to discuss the relationship between equiangular and biangular harmonic frames, and, ultimately, present the results needed for the next section, we must discuss their frame angle sets.  As mentioned, the frame angle set of a harmonic frame depends on the combinatorial relationship between the subselection $\mathcal S$ and its ambient group $\mathbb G$.  In particular, well-known~\cite{MR2246267} combinatorial objects known as {\bf difference sets} play an important role in the this discussion; however, it is convenient if we begin with a generalization, the {\bf bidifference sets}. 

\begin{definition}
Let $\mathcal S \subset \mathbb G$, where $\mathbb G$ is an additively written abelian group of order $n$ with identity $e$ and where $\mathcal S=\{g_1,...,g_m\}$ is a subset of $m$ elements.  We say that $\mathcal S$ is  an {\bf $(n,m,l, \lambda, \mu)$-bidifference set for $\mathbb G$ relative to $\mathcal A$} if 
$\mathcal A \subset \mathbb G$ is a subset of order $l$ with $e \in \mathcal A$, 
every non-identity element of $a\in A$ can be expressed as $a = g_s - g_t$ in exactly $\lambda$ ways and
every element $b\in \mathcal B=\mathbb G \backslash \mathcal A$ can be expressed as $b = g_s - g_t$ in exactly $\mu$ ways.

\end{definition}

In \cite{2016arXiv161003142C}, the authors considered a hiearchy of bidifference sets, atop which the {\it honest} difference sets are a special case of all of the others.  We define some of the more well-known~\cite{MR2246267, MR1440858, MR1277942} classes from this hierarchy.

\begin{definition}
Let $\mathcal S \subset \mathbb G$, where $\mathbb G$ is an additively written abelian group of order $n$ with identity $e$ and where $\mathcal S=\{g_1,...,g_m\}$ is a subset of $m$ elements. 
\begin{enumerate}
\item
We say that $\mathcal S$ is an {\bf $(n,m,\lambda)$-difference set} if it is an $(n,m,l, \lambda,\mu)$ for $\mathbb G$ relative to some subset $\mathcal A$ with $\lambda=\mu$.
\item
 Suppose $\mathbb H$ is a subgroup of $\mathbb G$ of order $l$.  We say that $\mathcal S$ is  an {\bf $(n,m,l, \lambda, \mu)$-divisible difference set for $\mathbb G$ relative to $\mathbb H$} if it is an $(n,m,l, \lambda, \mu)$-bidifference set for $\mathbb G$ relative to $\mathbb H$.  
We say that $\mathcal S$ is  an {\bf $(n,m,l,  \mu)$-relative difference set for $\mathbb G$ relative to $\mathbb H$} if it is an $(n,m,l,0,\mu)$-divisible difference set for $\mathbb G$ relative to $\mathbb H$. 
\item
We say that $\mathcal S$ is  an {\bf $(n,m,\lambda, \mu)$-partial difference set for $\mathbb G$} if 
$\mathcal S$ is an $(n,m,l, \lambda, \mu)$-bidifference set for $\mathbb G$ relative to $\mathcal S \cup \{ 0_G \}$,
where $l = |\mathcal S \cup \{ 0_G \}|$. 
\end{enumerate}
\end{definition}

Note that, trivially (or vacuously), every difference set may be viewed as both a divisible difference set and as a partial difference set.  Over the last decade, a characterization of equiangular harmonic frames in terms of difference sets has become  well-known~\cite{MR1984549, XiaZhouGiannakis2005, Fickus:2015aa}.  Motivated by this and the hierarchy outlined in the previous definition, the authors of~\cite{2016arXiv161003142C} studied the frame angle sets of harmonic frames generated by these other types of bidifference sets.  We summarize these results~\cite{MR1984549, XiaZhouGiannakis2005, 2016arXiv161003142C} in the following theorem.  

\begin{theorem}\label{thm_harmbtf}[\cite{MR1984549, XiaZhouGiannakis2005, 2016arXiv161003142C}]
Let $\mathbb G$ be an abelian group of order $n$ and let $\mathcal H$ be an $(n,m)$-harmonic frame generated by $\mathcal S$.
\begin{enumerate}
\item
If $\mathcal S$ is  an  $(n,m,l, \lambda, \mu)$-divisible difference set for $\mathbb G$ relative to $\mathbb H$, then $\mathcal H$ is a BTF with frame angle set
$$
\Theta_{\mathcal H} = \bigg\{\scriptstyle \frac 1 m \sqrt{m - \lambda + l(\lambda-\mu)}, \text{ } \frac 1 m \sqrt{m - \lambda }   \bigg\}.
$$
\item
If $\mathcal S$  is  an  $(n,m,\lambda, \mu)$-partial difference set for $\mathbb G$, then $\mathcal H$ is a BTF where the values of the frame angles are determined by the partial difference set's parameters (see~\cite{2016arXiv161003142C} for the frame angle formulae).
\end{enumerate}
Moreover, $\mathcal H$ is an ETF if and only if $\mathcal S$ is an $(n,m,\lambda)$-difference set for $\mathbb G$.
\end{theorem}

Many infinite families of difference sets~\cite{MR2246267}, divisible difference sets~\cite{ MR1440858}, and partial difference sets~\cite{MR1277942} are known within the combinatorial literature.  Accordingly, Theorem~\ref{thm_harmbtf} generates tables of infinite families of harmonic ETFs~\cite{Fickus:2015aa} and harmonic BTFs~\cite{2016arXiv161003142C}.  Although it is beyond the scope of this work to list all known examples here, we conclude this section by collecting three infinite families of harmonic BTFs, which we will use to construct {\it Steiner BTFs} in Section~\ref{sect_steiner}.

First, we note the existence of the underlying bidifference sets.

\begin{theorem}\label{thm_diff_set_exist}[\cite{MR2246267, MR1440858}; see also \cite{2015arXiv150905333B}]
\begin{enumerate}
\item {\bf Simplectic difference sets}
\\
 For every $n \in \mathbb N$ with $n>1$, an $(n, n-1, n-2)$-difference set exists.
\item {\bf Singer difference sets}
\\
For every prime power $q$, a $(q^2 + q + 1, q+1, 1)$-difference set exists.
\item  {\bf Picket fence sequences}
\\
For every prime power $q$, a $(q^2-1, q, q-1, 1)$-relative difference set for the additive group $\mathbb Z_{q^2-1}$ (relative to a subgroup $\mathbb H$ of order $q-1$) exists.
\end{enumerate}
\end{theorem}

Recalling that a relative difference set is a divisible difference set where the fourth parameter vanishes, ie $\lambda=0$, we apply Theorem~\ref{thm_harmbtf} to Theorem~\ref{thm_diff_set_exist} to obtain the desired families of BTFs.

\begin{corollary}\label{cor_harm_btfs_exists} $ $
\begin{enumerate}
\item {\bf Simplectic ETFs} \\
 For every $n \in \mathbb N$ with $n>1$, a complex, flat, equiangular $(n, n-1)$-frame $\mathcal H$ exists.
\item {\bf Singer ETFs} \\
For every prime power $q$ , a complex,
 flat, equiangular $(q^2 + q + 1, q+1)$-frame  $\mathcal H$ exists.
\item {\bf Picket fence BTFs} \\
 For every prime power $q$,  a complex, flat, biangular $(q^2-1, q)$-frame $\mathcal H$ exists with frame angle set
 $
\Theta_{\mathcal H} = \left\{\frac 1 q, \frac{1}{\sqrt q}  \right\}.
$
\end{enumerate}
\end{corollary}


\section{Steiner BTFs}\label{sect_steiner}
In this section, we construct numerous infinite families of biangular tight frames, called {\bf Steiner BTFs},  by exploiting the existence and ``flatness'' of harmonic BTFs along with well-studied combinatoral objects, called {\it Steiner  systems}~\cite{MR2246267}.
Our construction technique  is very similar to that of the {\bf Steiner ETFs} constructed in~\cite{MR2890902}, involving only a slight generalization of the so-called ``building blocks''. For those familiar with Steiner ETFs, we simply relax the requirement that the underlyinig difference sets are simplectic.

We remark that some of the families of Steiner BTFs in this section were previously described in a dissertation~\cite{haas_phd}, so we are happy to  present the results here in a more formal setting; however, we also augment these results~\cite{haas_phd} with new examples.  In particular, we introduce a second class of Steiner BTFs - distinct from those of \cite{haas_phd} -  which arise by passing from difference sets to certain types of divisible difference sets.


Roughly speaking, a {\it $(2,k,v)$-Steiner  system}, $(\mathcal V, \mathcal B)$, is a nonempty set, $\mathcal V$, of $v$ points along with a collection, $\mathcal B$, of $k$-subsets of $\mathcal V$ , called {\it blocks}, which satisfy certain incidence properties.  Every Steiner  system is associated to a $\{0,1\}$-matrix~\cite{MR2246267}, called its {\bf incidence matrix}.  For the purpose of constructing the {\it Steiner BTFs} in this section, we are mainly interested in the trasposes of the incidence matrices of such systems.  Accordingly, we find it less cumbersome  to simply assign a term and formal definition to the transpose of the incidence matrix of a Steiner  system - we call it a {\bf Steiner matrix} - with the tacit understanding that the existence of a Steiner matrix is equivalent to the existence of a Steiner  system~\cite{MR2246267}, where the parameters of the Steiner matrix completely determine those of the corresponding Steiner system.


\begin{definition}\label{def_steiner}
Suppose $k$ and $v$ are positive integers, where   $k\leq v$ and $\frac{v(v-1)}{k(k-1)}$ and $\frac{v-1}{k-1}$ are both integers.  A {\bf $(v,k)$-Steiner matrix} $A$  is a $\{0,1\}$-matrix  of size $\frac{v(v-1)}{k(k-1)} \times v $  such that
		\begin{enumerate}
		\item
		$A$ has exactly  $k$ ones in each row,
		\item
		$A$ has exactly  $\frac{v-1}{k-1}$ ones in each column, and
		\item 
		every two distinct columns of $A$ has a dot product of one.
\end{enumerate}
\end{definition}

As with the difference sets and their generalizations discussed in the previous section, constructions of numerous infinite families of $(2,k,v)$-Steiner  systems are known from the combinatorial literature~\cite{MR2246267}. Thus, the correspondence between Steiner systems and Steiner matrices outlined above implies the existence of numerous infinite families of corresponding Steiner matrices.  Their existence is vital to the main theorem of this section.

\begin{theorem}\label{thm_steiner_btfs}
		Let $A$ be a $(v,k)$-Steiner matrix.  To abbreviate, let $s=\frac{v-1}{k-1}$ and $m=\frac{v(v-1)}{k(k-1)}$, so that, by definition, $A$ is an $m \times v$ binary matrix  with exactly $s$ ones in each column. 
Let $H_1, H_2, ... , H_v$ 
be the synthesis matrices of flat $(t,s)$-frames for $\mathbb F^s$ such that either {\bf (i)} the corresponding frame $\mathcal H_j$ is equiangular  
for every $j \in \{1, 2, ..., v\}$ or  
{\bf (ii)}
 the corresponding frame $\mathcal H_j$ is biangular and $\frac 1 s \in \Theta_{\mathcal H_j}$  for every $j \in \{1, 2, ..., v\}$. 
Furthermore, let $F$ be the $m\times tv$ matrix constructed as follows:

		\begin{enumerate}
		\item
		For each $j \in \{1,.., v  \}$, let $F_j$ be the $m \times t$ (block) matrix obtained by replacing each one in the $j$th column of $A$ with a distinct row from $H_j$ and replacing each zero in the $j$th column of $A$ with a $1\times t$ row of zeros.
		\item
		Concatenate to obtain the $m \times tv$ matrix
		$F =[F_1 F_2 \cdots F_v]$.
		\end{enumerate}
		In this case, $F$ is the synthesis matrix of a biangular $(tv,m)$-frame, $\mathcal F$,  for $\mathbb F^m$.   
Moreover, the frame angle set of $\mathcal F$ is determined as follows:
\begin{enumerate}
\item[{\bf (i)}]
 If $\mathcal H$ is a flat ETF, then $\Theta_{\mathcal F} = \{\frac{1}{s}, W_{t,s} \}.$
\item[{\bf (ii)}]
 If $\mathcal H$ is a flat BTF with $\frac 1 s \in \Theta_{\mathcal H}$, then $\Theta_{\mathcal F} = \Theta_{\mathcal H}.$
\end{enumerate}

\end{theorem}
\begin{proof}
First, we show that $\mathcal F$ is indeed a $(tv, m)$-frame for  $\mathbb F^m$, noting that the following argument is independent of whether the $\mathcal H_j$s abides  Condition~{\bf (i)} or Condition~{\bf (ii)} in the hypothesis. 
Since every column of $F$ is really just a column of some $H_j$ inflated with extra zero entries, it follows that that $F$'s columns are unit-norm.
 It remains to verify that $\mathcal F$ is a tight frame for $\mathbb F^m$. Note that for any $j \in \{1,...,v\}$, the tightness of $\mathcal H_j$ implies that the inner product between any two distinct rows of $F_j$ must be zero.  Since the inner product between any two distinct rows of
$F$ is the sum of the inner products between the corresponding rows of the $F_j$s, it follows that 
the inner product between any two distinct rows of $F$ is zero.  Similarly, the squared norm of any row of $F$ is the sum of the corresponding squared norms of the rows of the $F_j$s, so the flatness of each $\mathcal H_j$ in conjunction with Property~(1) from Definition~\ref{def_steiner} implies that each row of $F$ has a squared norm of $\frac{kt}{s}$.  In particular, $FF^* = \frac{kt}{s} I_m$.  After noting the obvious dependence of the underlying field of $\mathcal F$ upon that of the  $\mathcal H_j$s and verifying the identity $\frac{kt}{s}=\frac{tv}{m}$, we conclude that $\mathcal F$ is a $(tv, m)$-frame for  $\mathbb F^m$.

	Next, we compute the frame angle set and verify that $F$ is a BTF.  Given two distinct columns $f$ and $f'$ of $F$,  there are two cases.  Either Case~(a)~they come from different block-matrices  (ie, $f$ is a column of $F_j$ and $f'$ is a column of $F_{j'}$ with $j \neq j'$) or Case~(b)~$f$ and $f'$ are columns of the same block-matrix $F_j$ for some $j$.  

 If Case~(a), then the flatness of $H_j$ and $H_{j'}$ along with Property~(3) from Definition~\ref{def_steiner} shows that $|\langle f, f' \rangle | = \frac{1}{s}$. In particular, $\frac 1 s \in \Theta_{\mathcal F}$.  Note that this is independent of whether the $\mathcal H_j$s satisfies  Condition~{\bf (i)} or Condition~{\bf (ii)}.

Note that either {\bf (i)} all of the $\mathcal H_j$s are flat ETFs or {\bf (ii)} they
are all flat BTFs with a shared frame angle of $\alpha_1 = \frac 1 s$, in which case the equidistributed property of BTFs along Equation~\ref{eq_unt_cond} implies that they must all agree on the second frame angle, $\alpha_2$.   In either case, we have
$$
\Theta_{\mathcal H_j} = \Theta_{\mathcal H_{j'}} \text{ for all } j,j' \in \{1,...,v\},
$$
so let us drop the unnecessary index $j$ and simply write $\Theta_{\mathcal H}=\Theta_{\mathcal H_j}$ for all $j \in \{1, 2, ..., v\}$.  

For Case~(b), it is clear by construction that $|\langle f, f' \rangle| \in \Theta_{\mathcal  H}$, so the arbitrariness of $j$ and $j'$ shows that $\Theta_{\mathcal F} = \left\{\frac 1 s\right\} \cup \Theta_{\mathcal H}$.
The claim follows by considering the two possible conditions from the hypothesis.
If ${\bf (i)}$ each $\mathcal H_j$ is an ETF, then $|\langle f, f' \rangle|=W_{t,s}$, the Welch bound.  Thus, $\Theta_{\mathcal F} = \{\frac{1}{s}, W_{t,s} \}$ and $\mathcal F$ is a BTF, as claimed.  On the other hand, if {\bf (ii)} each  $\mathcal H_j$ is biangular with $\frac 1 s \in \Theta_{\mathcal H_j}$, then $\Theta_{\mathcal F} = \left\{\frac 1 s\right\} \cup \Theta_{\mathcal H}=\Theta_{\mathcal H}$, and $\mathcal F$ is a BTF, as claimed.
\end{proof}

We call any BTF constructed by the preceding theorem a {\bf Steiner BTF}.  In order to demonstrate the theorem's potency, we construct Steiner BTFs by using harmonic ETFs and certain harmonic BTFs as  {\it building blocks}; more precisely, given a $(v,k)$-Steiner matrix $A$, we call a flat $(t,s)$-frame $\mathcal H$ a {\bf building block} for $A$ if $s=\frac{v-1}{k-1}$.
To begin, we recall the existence of four infinite families of $(2,k,v)$-Steiner systems and, hence, their corresponding $(v,k)$-Steiner matrices~\cite{MR2246267}.
\begin{theorem}\label{thm_steiner_sys}[\cite{MR2246267}] $ $
\begin{enumerate}
\item {\bf Affine geometries}
\\
Given a prime power $q$ and $a \in \mathbb N$ with $a\geq2$, then a $\left(q^a,q\right)$-Steiner matrix exists.
\item {\bf Projective geometries}
\\
Given a prime power $q$ and $a \in \mathbb N$ with $a\geq2$, then a $\left({\scriptstyle \frac{q^{a+1}-1}{q-1}, q+1}\right)$-Steiner matrix exists.
\item {\bf Unital systems}
\\
Given a prime power $q$, then a $\left(q^3+1, q+1\right)$-Steiner matrix exists.
\item {\bf Denniston systems}
\\
Given $a,b \in \mathbb N$ with $2 \leq a \leq b$, then
 a $\left(2^{a+b} + 2^a - 2^b, 2^a \right)$-Steiner matrix exists.
\end{enumerate}
\end{theorem}

In the following subsections, we use the harmonic ETFs and BTFs described in Section~\ref{sect_harmonic} as building blocks for these Steiner matrices, thereby producing numerous families of Steiner BTFs.

\subsection{Steiner ETFs}
Using the so-called {\it simplectic ETFs} from Corollary~\ref{cor_harm_btfs_exists} as building blocks in Theorem~\ref{thm_steiner_btfs}, we recover the plethora of {\bf Steiner ETFs} first described in~\cite{MR2890902}.

\begin{corollary}\label{cor_steiner_etfs}[\cite{MR2890902}]
If a $(v,k)$-Steiner matrix exists, then a complex, equiangular $(n,m)$-frame exists, where $m=\frac{v(v-1)}{k(k-1)}$ and $n=v(k+1)$.
\end{corollary}
\begin{proof}
Let $A$ be any $(v,k)$-Steiner matrix and, for brevity, write $s=\frac{v-1}{k-1}$. Let $\mathcal H$ be a complex, flat, equiangular $(s+1,s)$-frame, the existence of which is assured by Corollary~\ref{cor_harm_btfs_exists}.  Thus, $\mathcal H$ is an equiangular building block for $A$, so Theorem~\ref{thm_steiner_btfs} yields a Steiner BTF, $\mathcal F$, with angle set  $\Theta_{\mathcal F} = \{\frac{1}{s}, W_{s+1,s} \},$ so the claim follows by computing that $W_{s+1,s}= \frac 1 s$.
\end{proof}

Equiangular tight frames produced via Corollary~\ref{cor_steiner_etfs} are called {\bf Steiner ETFs}.  The authors of~\cite{MR2890902} provide a thorough analysis of the Steiner ETFs produced by each family of Steiner matrices from Theorem~\ref{thm_steiner_sys}.  We defer to their work~\cite{MR2890902} for further details.

\subsection{(Proper) Steiner BTFs from difference sets}
As noted in the previous section, a complete exhibition of all known harmonic ETFs is beyond the scope of this work; accordingly, we cannot provide a complete exhibition of all Steiner BTFs generated with (non-simplectic) harmonic ETFs as building blocks.  Instead,  we merely demonstrate the efficacy of Theorem~\ref{thm_steiner_btfs} by constructing three (of many possible) infinite families of Steiner BTFs via harmonic ETFs.

\begin{corollary}\label{cor_steiner_btfs_from_etfs}
If $q$ is a prime power, then
\begin{enumerate}
\item
{\bf Affine Singer (Steiner) BTFs} \\
a biangular $(n,m)$-frame $\mathcal F$ for $\mathbb C^{m}$ exists, where 
$m=q(q+1)$ and $n=q^2(q^2+q+1)$,
\item
{\bf Projective Singer (Steiner) BTFs} \\
a biangular $(n,m)$-frame $\mathcal F$ for $\mathbb C^{m}$ exists, where 
$m=q^2 +q + 1$ and $n=(q^2 +q + 1)^2=m^2$, and
\item
{\bf Unital Singer (Steiner) BTFs} \\
 if, in addition, $q$ is a square, then a
 biangular $(n,m)$-frame $\mathcal F$ for $\mathbb C^{m}$ exists, where 
$m=\frac{q^2(q^3+1)}{q+1}$ and $n=(q+1)(q^2+q+1)$.
\end{enumerate}

Moreover, in any of these three cases, the frame angle set of $\mathcal F$ is
$$\Theta_{\mathcal F} = \left\{ \frac{1}{q+1}, \frac{\sqrt q}{q+1} \right\}
.$$
\end{corollary}
\begin{proof}
Let $q$ be a prime power and let $\mathcal H$ be a complex, flat, equiangular $(q^2+q+1,q+1)$-frame, the existence of which is assured by Corollary~\ref{cor_harm_btfs_exists}.  By  Statement~(1) of Theorem~\ref{thm_steiner_sys}, a $(q^2, q)$ Steiner matrix $A$ exists, by Statement~(2), taking $a=2$, a $\left({ q^2+q+1, q+1}\right)$-Steiner matrix $B$ exists, and by Statement~(3), if $q$ is a square, then a $(q^3+1,q+1)$-Steiner matrix $C$ exists.  A straightforward computation shows that $\mathcal H$ is an equiangular building block for $A$, $B$, and $C$ (assuming $C$ exists).  Thus, Theorem~\ref{thm_steiner_btfs} yields a biangular $(n,m)$-frame $\mathcal F$ in each case, where the claimed values for $m$ and $n$  follow by elementary  computation.  The claimed value for the frame angle set in each case also follows from Theorem~\ref{thm_steiner_btfs}, after computing that $W_{q^2+q+1,q+1}= \frac{\sqrt q}{q+1}.$
\end{proof}

The Steiner BTFs  just described - produced with harmonic ETFs as building blocks - were previously mentioned in~\cite{haas_phd}.  We stress again that the BTFs produced in Corollary~\ref{cor_steiner_btfs_from_etfs} are far from comprehensive, as there are several other families of harmonic ETFs~\cite{MR2246267} and Steiner matrices~\cite{MR2246267} that can be used to produce more Steiner BTFs.
In the next subsection, we demonstrate a new class of Steiner BTFs, using certain  harmonic (proper) BTFs as building blocks instead.

\subsection{(Proper) Steiner BTFs from (proper) bidifference sets}
As with the previous section, a complete exhibition of all known Steiner BTFs generated by (proper) harmonic BTFs as building blocks is beyond the scope of this work.  In the following corollary, we use the so-called {\it picket fence sequences} of Corollary~\ref{cor_harm_btfs_exists} to demonstrate three families of Steiner BTFs produced via (proper) bidifference sets.

\begin{corollary}\label{cor_steiner_btfs_from_btfs}
If $q$ is a prime power such that $q+1$ is also a prime power, then
\begin{enumerate}
\item
{\bf Affine picket-fence (Steiner) BTFs} \\
a biangular $(n,m)$-frame $\mathcal F$ for $\mathbb C^{m}$ exists, where 
$m=q(q+1)$ and $n=q^3(q+1)$,
\item
{\bf Projective picket-fence (Steiner) BTFs} \\
a biangular $(n,m)$-frame $\mathcal F$ for $\mathbb C^{m}$ exists, where 
$m=q^2+q+1$ and $n=q(q+2)(q^2+q+1)$, and
\item 
{\bf Unital picket-fence (Steiner) BTFs} \\
if, in addition, $q$ is a square, then a
 biangular $(n,m)$-frame $\mathcal F$ for $\mathbb C^{m}$ exists, where 
$m=\frac{q^2(q^3+1)}{q+1}$ and $n=q(q+1)(q^3+1)$.
\end{enumerate}

Moreover, in any of these three cases, the frame angle set of $\mathcal F$ is
$$\Theta_{\mathcal F} = \left\{  \frac{1}{q+1}, \frac{1}{\sqrt{q+1}} \right\}
.$$
\end{corollary}

\begin{proof}
Let $q$ be a prime power such that $q+1$ is also a prime power, and let $\mathcal H$ be a complex, flat, biangular $\left((q+1)^2-1,q+1\right)$-frame with frame angle set $$\Theta_{\mathcal H} = \left\{ \frac{1}{q+1}, \frac{1}{\sqrt{q+1}} \right\},$$ the existence of which is assured by Corollary~\ref{cor_harm_btfs_exists}.  By  Statement~(1) of Theorem~\ref{thm_steiner_sys}, a $(q^2, q)$ Steiner matrix $A$ exists, by Statement~(2), taking $a=2$, a $\left({q^2 + q + 1, q+1}\right)$-Steiner matrix $B$ exists, and by Statement~(3), if $q$ is a square, then a $(q^3+1,q+1)$-Steiner matrix $C$ exists.  A straightforward computation shows that $\mathcal H$ is an biangular building block for $A$, $B$, and $C$ (assuming $C$ exists).  Thus, Theorem~\ref{thm_steiner_btfs} yields a biangular $(n,m)$-frame $\mathcal F$ with $\Theta_{\mathcal F}=\Theta_{\mathcal H}$ in each case, where the claimed values for $m$ and $n$ follow by elementary computation.
\end{proof}

Surprisingly, the preceding corollary along with an unintentional challenge implied by a work of Delsarte, Goethals, and Seidel~\cite{DelsarteGoethalsSeidel1975} leads to a connection with an open problem from number theory.

\subsubsection{Steiner BTFs, Mersenne primes, and Fermat primes}
We were careful not to refer to the families of Steiner BTFs produced by Corollary~\ref{cor_steiner_btfs_from_btfs} as {\it infinite families}, because
the necessary condition for such a frame's existence - the existence of a primer power $q$ such that $q+1$ is also a prime power - is equivalent to two famous open number-theoretic problems~\cite{MR990017}, namely the question of whether an infinite number of {\it Mersenne primes} or {\it Fermat primes} exist.  

To see this, note that if a pair, $q$ and $q+1$, of consecutive prime powers exists, then one of them must be even, and therefore a power of $2$.  In 2002, Mih\u{a}ilescu proved~\cite{MR2076124} the classical conjecture of Catalan~\cite{MR1578392}.

\begin{theorem}\label{thm_cat}[Catalan/Mih\u{a}ilescu; \cite{MR2076124, MR1578392}]
Given positive integers $a,b \in \mathbb N$ with $a,b\geq 2$, then the only solution to the Diophantine equation
$$
x^a +1 = y^b
$$
is $x=2,y=3,a=3$ and $b=2$.
\end{theorem}

Thus, if $q$ and $q+1$ are a pair of consecutive prime powers, then either Case~{\bf (i)} $q=8$ or Case~{\bf (ii)} one of the pair is an ``honest'' prime number while the other is a power of $2$.  As it turns out, Case~{\bf (ii)} boils down to $q$ either being a {\bf Mersenne prime}, a prime of the form $p=2^t-1$, or a  {\bf Fermat prime}, a prime of the form $p=2^t+1.$  Unfortunately, the questions of whether an infinitude of either type of primes exists have remained a famously open problems for centuries~\cite{MR990017}.  
Therefore, we cannot conclude that any of the three families of Steiner BTFs from Corollary~\ref{cor_steiner_btfs_from_btfs} are infinite; nevertheless, the corollary is obviously not vacuous, since, for example, the hypothesis is satsified for $q=2,3,4,7,8$ or any of the other known Merssenne primes~\cite{CW91} or Fermat primes~\cite{MR990017}.

To conclude this section, we briefly explain our ``challenge'' inferred from the work of Delsarte, Goethals, and Seidel~\cite{DelsarteGoethalsSeidel1975} which has motivated this number-theoretic digression.
In~\cite{DelsarteGoethalsSeidel1975}, the authors studied the upper bounds for the cardinalities of sets of unit vectors admitting few pairwise absolute inner products. Of relevance here,  they showed that for any complex $(n,m)$-frame $\mathcal F$, a necessary condition for $\mathcal F$ to be equiangular is that $n\leq m^2$ and a necessary condition for $\mathcal F$ to be biangular is that $n\leq \left(\tiny \begin{array}{c} m+1 \\ 2\end{array} \right)^2$ ; however,
besides a finite number of sporadic instances~\cite{Hoggar1982, 2016arXiv161003142C} and the known infinite families of {\it maximal sets of mutually unbiased bases} (see ~\cite{2016arXiv161003142C, WoottersFields1989} for details), we have found it strikingly difficult to construct or locate more complex BTFs in the literature~\cite{2016arXiv161003142C,Neumaier1989, Hoggar1982, WoottersFields1989, MR3047911, 2014arXiv1402.3521B} for which $n>m^2$; in other words, roughly speaking, it seems that most BTFs do not even exceed the cardinality bounds of ETFs.   Thus, we pose the following questions.

\begin{question}\label{quest_sq_btfs}
Other than maximal sets of mutually unbiased bases~\cite{WoottersFields1989}, do there exist other ``infinite families'' of complex biangular $(n,m)$-frames for which $n>m^2$? 
\end{question}

In light of the preceding discussion, if an infinite numer of Mersenne or Fermat primes exist, then the first and second families of Steiner BTFs produced by Corollary~\ref{cor_steiner_btfs_from_btfs} are infinite families; in particular, the biangular $\left({\scriptstyle q(q+2)(q^2+q+1), q^2+q+1}\right)$-frames produced by the second family (the projective picket fence BTFs) would answer Question~\ref{quest_sq_btfs} in the affirmative.


\section{Pl{\"u}cker ETFs}\label{sect_plucker}
In this final section,  we relax our emphasis from biangular tight frames to a well-studied generalization~\cite{2016arXiv160704546B, MR2520007, MR0313926, BachocEhler2013}, {\it (chordally) biangular tight fusion frames (BTFFs)}, with the goal of showcasing a surprising example of a BTFF which generates a Steiner ETF via the {\it Pl{\"u}cker embedding}.  To begin, we present a few bare essential facts about {\it fusion frames} and the {\it Pl{\"u}cker embedding.}

\subsection{Basics of fusion frames}
Let $l,m,n \in \mathbb N$, and 
let $\mathbb F = \mathbb R$ 
or $\mathbb F = \mathbb C$.
A {\bf (real or complex) tight $(n,l,m)$-fusion frame} is a set $\mathcal F=\{ P_j\}_{j=1}^n$, 
where each $P_j$ is an $m \times m$ orthogonal projection of rank $l$ with entries over $\mathbb F$ that satisfies the resolution of the identity,
\begin{equation}\label{eq_tff_id}
\sum_{j=1}^n P_j = a I_m, 
\end{equation}
for some  $a>0$, the fusion frame's  {\bf tightness parameter}.  
Recall that each $P_j$ satisfies $P=P^2=P^*$ and $\TR(P_j)=l$, so taking the trace of both sides of Equation~\ref{eq_tff_id} and solving yields the tightness parameter, $a=\frac{nl}{m}$.

Given a tight {\boldmath$(n,l,m)$}-fusion frame $\mathcal F$, its {\bf chordal fusion frame angles} are the elements of its {\bf fusion frame angle set},
$$
    \Theta_{\mathcal F} = \left\{ \sqrt{\TR\left(P_j P_{j'}\right) } : j \neq j' \right\}.
$$
We say that $\mathcal F$ is {\bf chordally $d$-angular} if $\left| \Theta_{\mathcal F} \right|=d$.  
In particular, if $d=1$ or $d=2$, 
then we call $\mathcal F$ a {\bf (chordally) equiangular tight fusion frame (ETFF)}
 or a {\bf (chordally) biangular tight fusion frame (BTFF)}, respectively.

Given any $(n,m)$-frame for $\mathbb F^m$, $\{f_j\}_{j=1}^n$, then it identifies with a tight $(n, 1, m)$-fusion frame for $\mathbb F^m$, $\{P_j\}_{j=1}^n$, where $P_j = f_j  f_j^*$ for each $j$. Moreover, by the identity
$$
\TR\left(P_j P_{j'}\right) = 
\TR\left( f_j f_j^*  f_{j'}  f_{j'}^*\right)
=
\left|\langle
f_j, f_{j'} \rangle
\right|^2
,$$
if $\{f_j\}_{j=1}^n$ is an ETF or BTF, then $\{P_j\}_{j=1}^n$ is ETFF or BTFF, respectively.  In this sense, one may regard the upcoming discussion of BTTFs as a natural generalization from that of  BTFs.

Finally, we recall a standard fact~\cite{} that identifies orthogonal projection matrices with equivalence classes of $1$-tight frames.

\begin{prop}\label{prop_projectionsandframes}[\cite{haas_phd, MR2964005}]
An $m \times m$ matrix $P$ is an orthogonal projection of rank $l$ over the field $\mathbb F$ if and only if there exists a synthesis matrix $F$ of a $1$-tight frame for $\mathbb F^l$ consisting of $m$ vectors such that $P=F^*F$.  Moreover, if $F'$ is another synthesis matrix of a $1$-tight frame for $\mathbb F^l$ consisting of $m$ vectors such that $P=(F')^*{F'}$, then there exists a unitary matrix, $U$, such that $F=UF'$.
\end{prop} 

In terms of the preceding proposition, we say that the frames $F$ and $F'$ are {\bf generators} of the projection $P$, and we say they are {\bf positively equivalent} if $\det(U)=1.$

\subsection{The Pl{\"u}cker embedding}
The {\it  Pl{\"u}cker embedding} is an algebrogeometric tool  used to study the Grassmannian manifold~\cite{MR507725, MR840877}, and it has found applications in frame theory several times~\cite{MR2964009, MR3085820}.  For the sake of brevity, we define it in terms of frame theory.

Let
 $\Omega_{m,l} (\mathbb F)$ denote the space of all $l \times m$ synthesis matrices of $1$-tight frames for $\mathbb F^l$ consisting of $m$ vectors.  Let $\rho_{m,l} = \left( \begin{array}{cc} m \\ l \end{array} \right)$, the number of ways to choose distinct $l\times l$ submatrices from an $l \times m$ matrix.  
After fixing an ordering,
 $$\{A_1, A_2, ..,A_{\rho_{m,l}}\},$$
on the $m\times m$ submatrices of an arbitrary $m \times n$ matrix $A$, we define the {\bf Pl{\"u}cker embedding}  as the coordinate mapping
$$
\Phi: \Omega_{m,l} \left( \mathbb F \right) \rightarrow \mathbb F^{\rho_{m,l}} : A \mapsto
\big[
\det (A_j)
\big]_{j=1}^{\rho_{m,l}}. 
$$
In light of Proposition~\ref{prop_projectionsandframes}, we can lift the  Pl{\"u}cker embedding's domain  to the space of all (real or complex) $m \times m$ orthogonal projections of rank $l$.  
Let $\mathcal G_{m,l} (\mathbb F)$ 
denote the space of all $m \times m$ orthogonal projections of rank $l$ over the field $\mathbb F$.  
For each   $P \in \mathcal G_{m,l} (\mathbb F)$,  
select a generator $A^{(P)} \in \Omega_{m,l} \left( \mathbb F \right)$ for $P$, and let $\left[\bm A^{(\bm P)} \bm\right]$ denote the class of all generators that are positively equivalent to $A^{(P)}$.
We define the {\bf lifted Pl{\"u}cker embedding} as the coordinate 
mapping
$$
\overline{\Phi} :  \mathcal G_{m,l} (\mathbb F) \rightarrow \mathbb F^{\rho_{m,l}}
 : P 
\mapsto
\Phi(A), \text{ for any } A \in  \left[\bm A^{(\bm P)} \bm\right]
$$
This is well-defined since, 
if $A^{(P)} \in  \left[\bm A^{(\bm P)} \bm\right]$
 and ${A'}^{(P)}=UA \in  \left[\bm A^{(\bm P)} \bm\right]$ 
are equivalent generators for $P$, then block identity,  
$$\det
 \left( {A'}^{(P)}_j \right) 
=\det(U) \det \left( A^{(P)}_j \right) 
= \det \left( A^{(P)}_j \right),
$$
shows the invariance of the value of $\overline{\Phi}(P)$ with respect to the choice of generator, $A \in   \left[\bm A^{(\bm P)} \bm\right]$.

\subsection{ A Pl{\"u}cker ETF}
We arrive at the main purpose of this section: an example of a BTFF that ``Pl{\"u}cker embeds'' into a Steiner ETF; in particular, we construct a chordally biangular (4, 2, 16)-tight fusion frame for $\mathbb R^4$, $\mathcal F = \{P_j\}_{j=1}^{16}$, such the  Pl{\"u}cker embedding of its elements,  $\overline{\mathcal F} = \big\{ \overline{\Phi}(P_j)\big\}_{j=1}^{16}$, forms a Steiner ETF for $\mathbb R^6$.
To begin, let 
$$
A =  \sqrt{\frac 2 3} \left[ 
\begin{array}{cccc}  
1 &  -1/2           &     -1/2            & 0 \\
0 &  \sqrt{3}/2   &   -\sqrt{3}/2     & 0
\end{array}
\right],
$$
and note that $A \in  \Omega_{4,2} (\mathbb R)$, since $FF^*= 2 I_2.$ Furthermore, let
$$
\Gamma = 
\left\{  c_j =
{\scriptscriptstyle
\begin{bmatrix}
\scriptscriptstyle 0 & \scriptscriptstyle 1 & \scriptscriptstyle 0 & \scriptscriptstyle 0 \\
\scriptscriptstyle 0 & \scriptscriptstyle 0 & \scriptscriptstyle 1 & \scriptscriptstyle 0 \\
\scriptscriptstyle 0 & \scriptscriptstyle 0 & \scriptscriptstyle 0 & \scriptscriptstyle 1 \\
\scriptscriptstyle 1 & \scriptscriptstyle 0 & \scriptscriptstyle 0 & \scriptscriptstyle 0
\end{bmatrix}
}^j :  j \in \{0, 1, 2, 3\}
\right\},
$$
the group of $4 \times 4$ cyclic matrix permutation matrices, and let 
$$
\Delta = 
\left\{  \scriptstyle
d_0 =
{
\scriptscriptstyle
	\begin{bmatrix}
	\scriptscriptstyle 1 & \scriptscriptstyle 0 & \scriptscriptstyle 0 & \scriptscriptstyle 0 \\
	\scriptscriptstyle 0 & \scriptscriptstyle 1 & \scriptscriptstyle 0 & \scriptscriptstyle 0 \\
	\scriptscriptstyle 0 & \scriptscriptstyle 0 & \scriptscriptstyle 1 & \scriptscriptstyle 0 \\
	\scriptscriptstyle 0 & \scriptscriptstyle 0 & \scriptscriptstyle 0 & \scriptscriptstyle 1
	\end{bmatrix}
},
d_1
=
{
\scriptscriptstyle
	\begin{bmatrix}
	\scriptscriptstyle 1 & \scriptscriptstyle 0 & \scriptscriptstyle 0 & \scriptscriptstyle 0 \\
	\scriptscriptstyle 0 & \scriptscriptstyle 1 & \scriptscriptstyle 0 & \scriptscriptstyle 0 \\
	\scriptscriptstyle 0 & \scriptscriptstyle 0 & \scriptscriptstyle -1 & \scriptscriptstyle 0 \\
	\scriptscriptstyle 0 & \scriptscriptstyle 0 & \scriptscriptstyle 0 & \scriptscriptstyle -1
	\end{bmatrix}
}, 
d_2=
{
\scriptscriptstyle
	\begin{bmatrix}
	\scriptscriptstyle 1 & \scriptscriptstyle 0 & \scriptscriptstyle 0 & \scriptscriptstyle 0 \\
	\scriptscriptstyle 0 & \scriptscriptstyle -1 & \scriptscriptstyle 0 & \scriptscriptstyle 0 \\
	\scriptscriptstyle 0 & \scriptscriptstyle 0 & \scriptscriptstyle 1 & \scriptscriptstyle 0 \\
	\scriptscriptstyle 0 & \scriptscriptstyle 0 & \scriptscriptstyle 0 & \scriptscriptstyle -1
	\end{bmatrix}
}, 
d_3=
{
\scriptscriptstyle
	\begin{bmatrix}
	\scriptscriptstyle 1 & \scriptscriptstyle 0 & \scriptscriptstyle 0 & \scriptscriptstyle 0 \\
	\scriptscriptstyle 0 & \scriptscriptstyle -1 & \scriptscriptstyle 0 & \scriptscriptstyle 0 \\
	\scriptscriptstyle 0 & \scriptscriptstyle 0 & \scriptscriptstyle -1 & \scriptscriptstyle 0 \\
	\scriptscriptstyle 0 & \scriptscriptstyle 0 & \scriptscriptstyle 0 & \scriptscriptstyle 1
	\end{bmatrix}
}
\right\}, 
$$
a diagonal unitary representation of the noncyclic group on four elements.  Next, let 
$$
\mathcal A =
\big\{
 A_{j,k} = A \, d_j \, c_k : j,k \in \{0,1,2,3\} 
\big\},
$$
so $\mathcal A$ is the orbit of $A$ under $\Gamma \times \Delta$.  Note that since the $c_j$s and $d_j$s are unitaries, it follow by the definition of $1$-tightness that $\mathcal A \subset \Omega_{4,2} (\mathbb R)$.  Thus, every element of $\mathcal A$ is the generator of some $4\times 4$ orthogonal projection of rank $2$.  Accordingly, we define the corresponding orthogonal projections,
$$
\mathcal P = \Big\{P_{j,k} = A_{j,k}^* A_{j,k} : j,k \in \{0,1,2,3 \} \Big\}.
$$
Next, we show that $\mathcal P$ is a chordally biangular tight $(16,2,4)$-fusion frame.  We have
$$
AA^* 
=
{
     \scriptscriptstyle
	\begin{bmatrix}
	\scriptscriptstyle 2/3 & \scriptscriptstyle -1/3 & \scriptscriptstyle -1/3 & \scriptscriptstyle 0 \\
	\scriptscriptstyle -1/3 & \scriptscriptstyle 2/3 & \scriptscriptstyle -1/3 & \scriptscriptstyle 0 \\
	\scriptscriptstyle -1/3 & \scriptscriptstyle -1/3 & \scriptscriptstyle 2/3 & \scriptscriptstyle 0 \\
	\scriptscriptstyle 0 & \scriptscriptstyle 0 & \scriptscriptstyle 0 & \scriptscriptstyle 0
	\end{bmatrix} 
},
\text{ and then, }
\sum\limits_{j}^4  d_j^* A^* A
 d_j
=
{
\scriptscriptstyle
	\begin{bmatrix}
	\scriptscriptstyle 8/3 & \scriptscriptstyle 0 & \scriptscriptstyle 0 & \scriptscriptstyle 0 \\
	\scriptscriptstyle 0 & \scriptscriptstyle 8/3 & \scriptscriptstyle 0 & \scriptscriptstyle 0 \\
	\scriptscriptstyle 0 & \scriptscriptstyle 0 & \scriptscriptstyle 8/3 & \scriptscriptstyle 0 \\
	\scriptscriptstyle 0 & \scriptscriptstyle 0 & \scriptscriptstyle 0 & \scriptscriptstyle 0
	\end{bmatrix}
},
$$
and then 
$$
\sum\limits_{j,k}^4
P_{j,k} =
\sum\limits_{j,k}^4 c_k^* d_j^* A^* A d_j c_k
=
\sum\limits_{k=1}^4 c_k^* {
\scriptscriptstyle
	\begin{bmatrix}
	\scriptscriptstyle 8/3 & \scriptscriptstyle 0 & \scriptscriptstyle 0 & \scriptscriptstyle 0 \\
	\scriptscriptstyle 0 & \scriptscriptstyle 8/3 & \scriptscriptstyle 0 & \scriptscriptstyle 0 \\
	\scriptscriptstyle 0 & \scriptscriptstyle 0 & \scriptscriptstyle 8/3 & \scriptscriptstyle 0 \\
	\scriptscriptstyle 0 & \scriptscriptstyle 0 & \scriptscriptstyle 0 & \scriptscriptstyle 0
	\end{bmatrix}
} c_k
=
8 \, I_4,
$$
 verifying that $\mathcal P$ is a real, tight $(16, 2, 4)$-fusion frame.  To see that $\mathcal P$ is chordally biangular, note that by our definition of the $P_{j,k}$s, its chordal frame angle set is
$$
\Theta_{\mathcal P} =
\Big\{
\TR
\left(
  c_k^* d_j^* A A^*  d_j d_ k       c_{k'}^* d_{j'}^* A A^*  d_{j'} d_{k'}\right) : j, j', k, k' \in \{0, 1, 2, 3\}
\Big\}.
$$
Given $j, j', k, k' \in \{0, 1, 2, 3\},$
  if $k=k'$ but $j \neq j'$, then a straighforward computation yields
$$
\TR (P_{j, k} P_{j', k}) = \frac{10}{9}.
$$
Otherwise, if $j \neq j'$ and $k \neq k'$, 
the the computation reduces to the trace inner product between two $2\times 2$ principle submatrices,
$$
\TR \left(P_{j, k} P_{j', k'} \right) = 
\TR
\left(
	\left[
		\begin{array}{cc} 
		 2/3 & \pm 1/3 \\ \pm 1/3 & 2/3
		\end{array} 
	\right]
	\left[
		\begin{array}{cc} 
		 2/3  & \pm 1/3 \\ \pm 1/3 & 2/3
		\end{array} 
	\right]
\right)
=
\frac 4 9 \pm \frac 2 9.
$$
Thus, $\mathcal P$ is a BTFF with chordal frame angle set
$$
\Theta_{\mathcal P} = \left\{ \sqrt{10/9}, \sqrt{2/3}   \right\}.
$$

Finally, we show that $\mathcal P$  ``Pl{\"u}cker  embeds'' into a Steiner ETF.  Viewing $\mathcal A$ as a set of generators for the elements of $\mathcal P$, we compute the Pl{\"u}cker  embeddings, ${\mathcal F} = \left\{  \overline{\Phi} (P_j) \right\}$, and conclude that $\mathcal F$ is a Steiner ETF.
Note that $\rho_{4,2}=6$, and recall that the Pl{\"u}cker embedding requires a choice of ordering for the underlying submatrices.  We use a basic dictionary ordering; that is, given an abitrary $2 \times 4$ matrix $A=[a_1 \, \, a_2 \, \, a_3  \, \, a_4]$ with columns $a_1, a_2, a_3$ and $a_4$, we order the $2 \times 2$ submatrices as 
$$\scriptstyle
A_1=[a_1  \, a_2], A_2=[a_1  \,  a_3], A_3=[a_1  \,  a_4], A_4=[a_2  \,  a_3], A_5=[a_2   \, a_4], \text{ and } A_6=[a_3   \, a_4].
$$
Thus,  computing the Plucker embedding of $P_{j,k}$ is a simple matter of computing the determinates of its six $2 \times 2$ submatrices of its generator, $A_{j,k}$.  Fixing $k=0$, define the $6 \times 4$ (block) matrix,
\begin{align*}
F_0 
&= \left[ \overline{\Phi} (P_{j,0}) \right]_{j=1}^6 \\
&= 
{\scriptscriptstyle\tiny
\sqrt{2} \left[
		\begin{array}{cccc} 
		       1 & 1  & -1 & -1 \\
                 -1 & 1 & -1 & 1 \\
                 0 &  0  & 0 & 0 \\
                 1 & -1 & -1 & 1 \\
                 0 &  0  & 0 & 0 \\
                 0 &  0  & 0 & 0 \\
		\end{array} 
\right]}
\end{align*}
Similarly, define and compute $F_k = \left[ \overline{\Phi} (P_{j,k}) \right]_{j=1}^6$ for $k=1,2$ and $3$, 
$$\tiny\scriptstyle
F_2
=
{\scriptscriptstyle\tiny
\sqrt{2} \left[
		\begin{array}{cccc} 
		       0 &  0  & 0 & 0 \\
                 0 &  0  & 0 & 0 \\
                 0 &  0  & 0 & 0 \\
                  1 & 1  & -1 & -1 \\
                 -1 & 1 & -1 & 1 \\
                 1 & -1 & -1 & 1 \\
		\end{array} 
\right]
},
F_3
=
{\scriptscriptstyle\tiny
\sqrt{2} \left[
		\begin{array}{cccc} 
                  0 &  0  & 0 & 0 \\
		       1 & -1  & 1 & -1 \\
                 -1 & 1 & 1 & -1 \\
                 0 &  0  & 0 & 0 \\
                 0 &  0  & 0 & 0 \\
                 1 & 1 & -1 & -1
           \end{array} 
\right]
},
\text{ and }
F_4
=
{\scriptscriptstyle\tiny
\sqrt{2} \left[
		\begin{array}{cccc} 
		       1 & -1  & -1 & 1 \\
                 0 &  0  & 0 & 0 \\
                 -1 & -1 & 1 & 1 \\
                 0 &  0  & 0 & 0 \\
                 1 & -1 & 1 & -1 \\
                 0 &  0  & 0 & 0 \\
		\end{array} 
\right]
}.
$$
Concatenate and rescale to obtain the $6 \times 16$ matrix
$$
F = \frac{1}{\sqrt 6} \Big[ F_1 \, F_2 \, F_3 \, F_4 \Big].
$$

Notice that for each $j\in \{0,1,2,3\}$, the $3\times 4$ submatrix, $H_j$, obtained by deleting the three zero rows from $F_j$ and rescaling appropriately, forms the synthesis matrix of a real, flat equiangular $(4,3)$-frame.
By inspection, it is clear that $F$ is constructed in accordance with Theorem~\ref{thm_steiner_btfs} with respect to the 
$(6,2)$-Steiner matrix
$$S= \tiny \scriptscriptstyle
\begin{bmatrix}
1 & 0 & 0 & 1 \\
1 & 0 & 1 & 0 \\
0 & 0 & 1 & 1 \\
1 & 1 & 0 & 0 \\
0 & 1 &  0 & 1 \\
0 & 1 & 1 & 0
\end{bmatrix},
$$
using the $H_js$ as building blocks.  Noting that $W_{4,3}=\frac 1 3$, we conclude that $F$ is the synthesis matrix of a Steiner ETF, $\mathcal F$, for $\mathbb R^6$ consisting of $16$ vectors.
In light of this construction, we conclude with a natural question, the posing of which requires a definition.
\begin{definition}
An equiangular $(n,m)$-frame for $\mathbb F^m$, $\mathcal H$, with synthesis matrix, $H=[h_1 \, h_2 \, ... \, h_n]$, is a {\bf Pl{\"u}cker ETF} if there exists a (real or complex) tight $(n, l, t)$-fusion frame $\mathcal Q = \{Q_j\}_{j=1}^n$, where $l\geq2$, $m=\rho_{t,l}$ and $f_j = \overline{\Phi} \left(Q_j\right)$ for every $j\in\{1,2,...,n\}$.  
\end{definition}

Thus, the Steiner ETF, $\mathcal F$, that we have just constructed is a Pl{\"u}cker ETF.  Given the volumetric nature of the Pl{\"u}cker  embedding (ie, the determinants of the submatrices correspond to (signed) hyper-volumes of hyper-parallopipeds) and that - as an astute reader might have noticed - the generators for each $P_j$ from our example are essentially $2$-simplices living in $\mathbb R^4$, we find our construction strangely intuitive.  Nevertheless, we are unaware of any other examples.
\begin{question}
Other that the Pl{\"u}cker ETF, $\mathcal F$, constructed in this section, do  Pl{\"u}cker ETFs with different parameters exist?  
\end{question}
%


\bibliographystyle{amsplain}
\bibliographystyle{amsplain}
\bibliography{etfsvsbtfs}
\end{document}